\documentclass[12pt,reqno]{amsart}

\usepackage{pb-diagram,verbatim}

\newtheorem{theorem}{Theorem}[section]
\newtheorem*{thm5.1}{Theorem 5.1}
\newtheorem{lemma}[theorem]{Lemma}
\newtheorem{proposition}[theorem]{Proposition}
\newtheorem{corollary}[theorem]{Corollary}

\theoremstyle{definition}
\newtheorem{definition}[theorem]{Definition}
\newtheorem{remark}[theorem]{Remark}
\newtheorem{example}[theorem]{Example}

\numberwithin{equation}{section}

\author{Geoffrey Buhl}
\address{Department of Mathematics, California State University  Channel Islands, Camarillo, CA
93012}
\email{geoffrey.buhl@csuci.edu}

\author{Gizem Karaali}
\address{Department of Mathematics, Pomona College, Claremont, CA
91711}
\email{gizem.karaali@pomona.edu}

\thanks{Geoffrey Buhl is supported by a California State University Channel Islands Faculty Development Grant.}


\date{\today.}


\theoremstyle{definition}

\newcommand{\baseRing}[1]{\ensuremath{\mathbb{#1}}}

\newcommand{\Z}{\baseRing{Z}}
\newcommand{\N}{\baseRing{N}}

\renewcommand{\d}{\delta}

\renewcommand{\phi}{\varphi}
\newcommand{\Res}{\operatorname{Res}}
\newcommand{\Span}{\operatorname{span}}
\newcommand{\wt}{\operatorname{wt}}
\newcommand{\row}{\operatorname{row}}
\newcommand{\End}{\operatorname{End}}\newcommand{\1}{{\textbf 1}}

\allowdisplaybreaks

\begin{document}

\title[Spanning sets for M\"{o}bius vertex algebras]{Spanning sets for M\"{o}bius vertex algebras satisfying arbitrary difference conditions}

\begin{abstract}
Spanning sets for vertex operator algebras satisfying difference-zero and difference-one conditions have been extensively studied in the recent years. In this paper, we extend these results. More specifically, we show that for a suitably chosen generating set, any $\N$-graded M\"{o}bius vertex algebra is spanned by monomials satisfying a difference-$N$ ordering condition.

\end{abstract}
\maketitle

\section{Introduction}

The theory of vertex algebras and vertex operator algebras shares fundamental connections with number theory, the theory of finite simple groups, and string and conformal field theories in physics. These connections are manifested in the Moonshine Module vertex operator algebra, whose symmetry group is the Monster simple group and whose graded trace is the modular function $j(\tau)$.

In Zhu's work on the modularity properties of certain vertex operator algebras, a key assumption is \emph{$C_2$-cofiniteness} (also called \emph{Zhu's finiteness condition}): a special subspace of the vertex operator algebra, henceforth called {\it Zhu's subspace}, has finite codimension \cite{MR1317233}.  In \cite{MR1990879}, Gaberdiel and Neitzke use representatives of a quotient of a vertex operator algebra by Zhu's subspace to generate a Poincar\'{e}-Birkhoff-Witt-like spanning set for vertex operator algebras satisfying certain ordering restrictions. 

The existence of such a spanning set has various interesting implications for the representation theory of the associated vertex (operator) algebra. For instance, Gaberdiel and Neitzke use their spanning set to prove Nahm's conjecture: if Zhu's subspace has finite codimension then a certain subspace of any irreducible module has finite codimension \cite{MR1990879}. A generalization of the same spanning set for modules can be used to show that under the assumption that Zhu's subspace has finite codimension, two different notions of complete reducibility for vertex operator algebras are equivalent  \cite{MR2052955}.

In this paper, we generalize such spanning set results. In particular we show the following: Using representatives of a suitable generalization of Zhu's subspace, we can construct a family of Poincar\'{e}-Birkhoff-Witt-like spanning sets that satisfy certain ordering restrictions.  We construct these spanning sets for M\"{o}bius vertex algebras, which are certain generalizations of vertex operator algebras.  Rather than admiting a represention of the full Virasoro algebra, a M\"{o}bius vertex algebra admits a represention of $\mathfrak{sl}(2)$, a subalgebra of the Virasoro algebra.

We now go into some more technical details. 
In \cite{MR1700507} Karel and Li show that the set of representatives of a certain subspace of a vertex operator algebra is a minimal generating set for a spanning set for the vertex operator algebra. This spanning set consists of monomials in the generators which satisfy ordering restrictions similar to the Poincar\'{e}-Birkhoff-Witt bases of Lie algebras. In particular the relevant monomials are of the form
\begin{eqnarray}
x^{1}_{n_1}\cdots x^{k}_{n_k}\1,
\end{eqnarray}
where the $x$'s are elements of the generating set and $$\deg(x^{1}_{n_1})\geq \deg(x^{2}_{n_2}) \geq \cdots \geq \deg(x^{k}_{n_k}) > 0.$$  This ordering restriction is equivalent to $n_1\leq n_2 \leq \cdots \leq n_k < 0$ after suitable rearrangement.

In \cite{MR1990879} Gaberdiel and Neitzke construct a similar Poincar\'{e}-Birkhoff-Witt-like spanning set for a given vertex operator algebra with a larger set of generators. Their generators are representatives of a quotient space of the relevant vertex operator algebra by Zhu's subspace. Basically, their spanning set consist of monomials of the form
\begin{eqnarray}
x^{1}_{n_1}\cdots x^{k}_{n_k}\1,
\end{eqnarray}
where the $x$'s are elements of a generating set larger than that of Karel and Li in \cite{MR1700507}. The indices $n_i$ of the modes have to satisfy a no-repetition restriction: $n_1< n_2 < \cdots < n_k < 0$. In other words, they must be strictly increasing. In short, by expanding the size of the generating set for this spanning set, Gaberdiel and Neitzke were able to allow for more restrictive order conditions on the monomials in the spanning set. 
This construction of the no-repetition spanning set has later on been extended to a spanning set for modules for vertex operator algebras \cite{MR1927435}, twisted modules for vertex operator algebras \cite{MR2039213}, and quasi-modules for M\"{o}bius vertex algebras \cite{quasi}.

We can reformulate order conditions on the modes of indices like the no-repetition restriction in a more natural way, using difference conditions.  In this framework, the Karel and Li spanning set obeys a difference-zero condition, and the Gaberdiel and Neitzke spanning set obeys a difference-one condition.
This perspective is motivated by the work of Lepowsky and Wilson on the Lie-theoretic interpretations of certain Rogers-Ramanjuan identities  \cite{MR638674, MR781375, MR752821}. 
There, they develop a ``straightening'' procedure to, in one instance, construct a basis for certain $A^{(1)}_1$-modules  that satisfy a difference-two condition.  From this perspective, we see the identity of Gaberdiel and Neitzke 
\[ u_{n}v_{n} = (u_{-1}v)_{2n+1} - \sum_{i \ge 0, i \neq -n} u_{-i} v_{2n+i} 
- \sum_{i \ge 0} v_{2n-i} u_i. \]
as a difference-one straightening identity for vertex algebras.  

In this paper we generalize the difference-one straightening identity above to a family of straightening identities. We then use these identities to build on the earlier spanning set results. We explicitly show that by increasing the size of the generating set, we can construct spanning sets for M\"{o}bius vertex algebras with more and more restrictive ordering conditions. 

Thus the main result of this paper can be summarized as follows: For suitably chosen generating sets, there are Poincar\'{e}-Birkhoff-Witt-like spanning sets for M\"{o}bius vertex algebras that satisfy a difference-$N$ condition.  More specifically, for a given M\"{o}bius vertex algebra $V$, there are spanning sets that consist of monomials of the form 
\begin{eqnarray}
x^{1}_{n_1}\cdots x^{k}_{n_k}\1,
\end{eqnarray}
where the $x$'s are elements of a specific generating set, and the $n_i$'s are negative integers with $n_{i+1} -n_i   \ge N$.  The generators of this difference-$N$ spanning set are representatives of a basis of a quotient $V/C_{N+1}(V)$ of $V$ by a generalization $C_{N+1}(V)=\{u_{-N-1}v: u,v \in V \}$ of Zhu's subspace. This main result of the paper is recorded as:

\begin{thm5.1} 
For an $\N$-graded M\"{o}bius vertex algebra $V$ and $N\in\Z_+$, let $X$ be a set of homogeneous representatives of a spanning set for the quotient space $V/C_{N+1}(V)$, $V$ is spanned by the elements of the form
\begin{eqnarray}
x^1_{n_1} x^2_{n_2} \cdots x^k_{n_k}\1
\end{eqnarray}
where $k \in \N$; $x^{1}, \ldots, x^{k} \in X$; $n_1, \ldots , n_k \in \Z_{-}$; and $n_1 < n_2 < \cdots < n_k < 0$ with $n_i - n_{i+1} \ge N$ for each $1\ge i \ge k-1$.
\end{thm5.1}

\noindent
In other words, $V$ is spanned by monomials satisfying the following order restriction: the indices of adjacent modes differ by at least $N$.

The structure of this paper is as follows:  In Section 2 we present some preliminary information including the definition of a M\"{o}bius vertex algebra. In Section 3 we define a family of filtrations for M\"{o}bius vertex algebras and explore some of the properties of these filtrations.  In Section 4 we prove a family of straightening lemmas used to prove our main theorem. Section 5 contains the statement and the proof of the main theorem.  Section 6, the appendix, contains some combinatorial arguments used in the construction of the straightening lemmas of section four.

\section{M\"{o}bius vertex algebras}

In this section we present the definitions of a M\"{o}bius vertex algebra, the $N$-th Zhu subspace, and other important concepts.  For a primer of the theory of vertex algebras and vertex operator algebras, we refer the reader to \cite{MR2023933}. 

The notion of a {\it M\"{o}bius vertex algebra} is a generalization of the notion of a {\it conformal vertex algebra}. The latter is a $\Z$-graded vertex algebra that admits a Virasoro algebra representation, while the former admits a representation of $\mathfrak{sl}(2)$, as a Lie subalgebra of the Virasoro algebra. M\"{o}bius vertex algebras appear in the work of Huang, Lepowsky, and Zhang on a logarithmic tensor product theory for conformal vertex algebras \cite{huang-2006}.  

\begin{definition}
A {\it M\"{o}bius vertex algebra} is a $\Z$-graded vector space
 \begin{eqnarray*}
V=\coprod_{n \in \Z} V_n
\end{eqnarray*}
equipped with a linear map 
\begin{eqnarray*}
Y:  V &\rightarrow& \End(V)[[x,x^{-1}]]\\
  v &\mapsto& Y(v,x)=\sum_{n \in \Z}v_n x^{-n-1} \ \mbox{(where $v_n \in \End(V)$)}
\end{eqnarray*}
where $Y(v,x)$ is called the {\it vertex operator associated with $v$} and a distinguished vector $\1 \in V_0$ (the {\it vacuum vector}), satisfying the following conditions for $u,v \in V$: 

\noindent
The lower truncation condition: 
\begin{eqnarray*}
u_nv=0 \mbox{ for n sufficiently large};
\end{eqnarray*}
the vacuum property:
\begin{eqnarray*}
Y(\1,x)=1_V;
\end{eqnarray*}
the creation property:
\begin{equation}
\label{CreationAxiom}
Y(v,x) \1 \in V[[x]] \mbox{ and } \lim_{x \rightarrow 0} Y(v,x) \1=v;
\end{equation}
and the Jacobi Identity:
\begin{eqnarray}
&x_0^{-1}\d \left(\frac{x_1 - x_2}{x_0}\right)Y(u,x_1)Y(v,x_2)-
x_0^{-1} \d \left(\frac{x_2- x_1}{-x_0}\right)Y(v,x_2)Y(u,x_1) \nonumber \\
\label{JacobiIdentity}
&=x_2^{-1} \d \left(\frac{x_1- x_0}{
x_2}\right)Y(Y(u,x_0)v,x_2).
\end{eqnarray}
In addition there is a representation $\rho$ of $\mathfrak{sl}(2)$ on V given by:
\begin{eqnarray*}
L(j)=\rho(L_j) \mbox{, } j=-1,0,1
\end{eqnarray*}
where $\{ L_{-1}, L_0, L_1\}$ from a basis of $\mathfrak{sl}(2)$ with Lie brackets
\begin{eqnarray*}
[L_0,L_{-1}]= L_{-1}, [L_0,L_1]=-L_1  \mbox{, and} [L_{-1},L_{1}]=-2L_0, 
\end{eqnarray*}
and the following conditions hold for $v \in V$ and $j=-1,0,1$:
\begin{eqnarray}
[L(j), Y(u,x)]&=& \sum^{j+1}_{k=0} \binom{j+1}{k}x^{j+1-k}Y(L(k-1)v,x), \nonumber\\
\frac{d}{dx}Y(v,x)&=&Y(L(-1)v,x), \label{l-1}
\label{LDerivationProperty}
\end{eqnarray}
and
\begin{eqnarray*}
L(0)v=nv=(\wt v)v \mbox{ for } n \in \Z \mbox{ and } v \in V_n.
\end{eqnarray*}
\end{definition}
\noindent
A M\"{o}bius vertex algebra is denoted by the quadruple $(V,Y,\1,\rho)$ or by $V$ when clear from the context. 

M\"{o}bius vertex algebras are  generalizations of {\it quasi-vertex operator algebras}. Quasi-vertex operator algebras  which are generalizations of vertex operator algebras, are defined in \cite{MR1142494}, and include two axioms in addition to those of a M\"{o}bius vertex algebra: each graded piece is finite dimensional, and the $\Z$-grading is truncated from below. Thus we would like to be able to talk about the M\"{o}bius vertex algebras which satisfy a similar truncation condition:

\begin{definition}
A M\"{o}bius vertex algebra $V=\coprod_{n \in \Z} V_n$ is {\it $\N$-graded} if $V_n=0$ for $n<0$.
\end{definition}

In the literature on vertex (operator) algebras, depending on the context of a problem, the type of vertex algebra studied may have: finite dimensional graded pieces, as in the case of quasi-vertex operator algebras and vertex operator algebras; a representation of the full Virasoro algebra, as in the case of conformal vertex algebras and vertex operator algebras; a lower truncation condition, as in the case of vertex operator algebras and quasi-vertex operator algebras; or an $\mathfrak{sl}(2)$ representation, as in the case of M\"{o}bius vertex algebras and quasi-vertex operator algebras.  In this paper, we require that the vertex algebras we study be $\N$-graded, satisfy a specific lower-truncation condition and have the $L(-1)$-derivation property (i.e. Equation \eqref{LDerivationProperty});  $\N$-graded M\"{o}bius vertex algebras fulfill these requirements. 

We now focus on certain subspaces of our vertex algebras:

\begin{definition}
For a M\"{o}bius vertex algebra $V$ and $N \in \Z_+$, the \emph{$N$-th Zhu subspace} is
\begin{eqnarray}
C_N(V)=\{u_{-N}v: u,v \in V \}.
\end{eqnarray}
\end{definition}

Note that for $N=1$, this definition differs from the $C_1$-subspace as defined in \cite{MR1700507} \cite{quasi}; in our case, the creation axiom  implies that $V=C_1(V)$.

The above generalizes the notion of {\it Zhu's subspace} $C_2(V)$.  Assuming the cofiniteness of this subspace is crucial to proving certain modularity properties of modules and twisted modules for vertex operator algebras and vertex operator superalgebras \cite{MR1317233,MR1794264,dong-2006,MR2175996}.  Cofiniteness of Zhu's subspace is called \emph{Zhu's finiteness condition} or \emph{$C_2$-cofiniteness} and is generalized in the following definition.

\begin{definition}
For  $N \geq 0$,  a vertex algebra $V$ is {\it $C_N$-cofinite}  if the $N$-th Zhu subspace $C_N(V)$ has finite codimension in $V$.
\end{definition}

Comparing coefficients of $x^{-N-1}$ in Equation (\ref{l-1}), we get 
\[-Nu_{N-1}=(L(-1)u)_N\] 
\noindent
which implies that $C_{N-1}(V) \subset C_N(V)$. So we have 
\[V=C_1(V) \supset C_2(V) \supset \cdots \supset C_N(V) \cdots\]
\noindent
If $V/C_{N}(V)$ is finite-dimensional then $V/C_{N-1}(V)$ is finite-dimensional.  In other words, $C_{N}$-cofiniteness implies $C_{N-1}$-cofiniteness. In \cite{MR1990879}, the difference-one spanning set for algebras is used to prove the converse: $C_{N-1}$-cofiniteness implies $C_{N}$-cofiniteness for $N \ge 3$.  

In this paper, we do not require $C_N$-cofiniteness (or equivalently, $C_2$-cofiniteness). However, if a M\"{o}bius vertex algebra is $C_2$-cofinite, then the spanning sets given by our main (difference-$N$ spanning set) theorem will be finitely generated.  In particular, see Corollary \ref{diffncor}.

We now look at some instances of the Jacobi identity (Equation \eqref{JacobiIdentity}). We will use two of its specializations in this paper: these will be the associativity and commutativity identities for modes of M\"{o}bius vertex algebras. These identities are also called \emph{Borcherds's iterate} and \emph{commutator formulas}, respectively \cite{MR843307}. 

{\it Associativity identity}: For  $u,v \in V$ a vertex algebra and $m,n \in \Z$,
\begin{equation}
\label{lemmaassociator}
(u_mv)_n = \sum_{i \ge 0} \begin{pmatrix} m \\ i \end{pmatrix} (-1)^i \left(u_{m-i} v_{n+i} -  (-1)^{m}v_{m+n-i} u_i \right).
\end{equation}
It is obtained by taking $\Res_{x_0} \Res_{x_1} \Res_{x_2} x_0^m x_2^n$ of the Jacobi identity.

{\it Commutativity identity}: For  $u,v \in V$ a vertex algebra and $m,n \in \Z$,
\begin{equation}
[u_m, v_n] = \sum_{j \ge 0} \begin{pmatrix}m\\j\end{pmatrix} (u_jv)_{m+n-j}. \label{lemmacommutator}
\end{equation}
It is obtained by taking $\Res_{x_0} \Res_{x_1} \Res_{x_2} x_1^m x_2^n$ of the Jacobi identity.

The associativity and commutativity identities for modes play a fundamental role in determining the properties of the difference-N filtrations defined in the next section.  The associativity identity is also used to construct the straightening lemmas of section four.

\section{Filtrations for M\"{o}bius vertex algebras}

A key ingredient in the proof of the difference-zero and difference-one spanning sets is a certain filtration of the vertex operator algebras.  The filtration is due to Watts an§d allows for certain rearrangement and replacement properties \cite{MR1068697}. In this section we generalize this filtration to obtain a family of increasing filtrations for graded vertex algebras. Each filtration in this family is good in the sense of Li \cite{MR2048777} and allows for some general rearrangement and replacement properties.

The filtration for a vertex operator algebra $V$ given by Watts is:
\[ V^{(0)} \subset V^{(1)} \subset V^{(2)} \subset \cdots \subset V, \]
where 
\[ V^{(s)} = \Span\{ u_{n_1}^1 u_{n_2}^2 \cdots u_{n_k}^k  1 : \sum_{i = 1}^k \wt{u^i}   \le s\}. \]
In our family of filtrations, this will be the first member, and will correspond to $N=1$.  The {\it rearrangement property} for this filtration is: given a monomial of filtration level $s$, any permutation of the modes in the monomial is equivalent to the original modulo a lower filtration element.  The {\it replacement property} for this filtration is: replacing a mode in a monomial of filtration level $s$ with another mode in the same coset of $V/C_2(V)$ results in an equivalent monomial modulo a lower filtration element.

\subsection{Difference-$N$ Filtrations}
  
For the difference-$N$ spanning set, we generalize the above filtration in the following way: the $s$-th filtered level contains all monomials of the form $u_{n_1}^1 u_{n_2}^2 \cdots u_{n_k}^k  1$ such that $ \sum_{i = 1}^k (\wt{u^i} + (N-1)) \le s$.

\begin{definition}
The {\it difference-$N$ filtration} of a M\"{o}bius vertex algebra $V$ is 
\[ V^{(0)}_N \subset V^{(1)}_N \subset V^{(2)}_N \subset \cdots \subset V, \]
where 
\[ V^{(s)}_N = \Span\{ u_{n_1}^1 u_{n_2}^2 \cdots u_{n_k}^k  1 :  (N-1)k + \sum_{i = 1}^k \wt{u^i}   \le s\}. \]
\end{definition}

This is a family of increasing filtrations on $V$, with the $M$-th filtration finer than the $N$-th filtration for $M>N$. For any fixed filtration level $s$, as $N$ increases, the size of $V^{(s)}_N$ decreases. That is,
\[ V^{(s)}_1 \supset V^{(s)}_2 \supset V^{(s)}_3 \supset \cdots \supset V^{(s)}_N \supset \cdots. \] 

\subsection{Properties of the Filtration}

The proof of the difference-one spanning set for vertex operator algebra relies on rearrangement and replacement properties of the difference-one filtration.  We now demonstrate that these properties hold for each difference-$N$ filtration.  It is these properties that provide the framework for the induction argument that is used to prove the difference-$N$ spanning set theorem.

The following lemma demonstrates that difference-$N$ filtrations satisfy the rearrangement property of the Watts filtration: rearranging the modes in a monomial does not introduce any terms of a higher filtration level. In other words, such rearrangements result in an equivalent monomial modulo a lower filtration level.

\begin{lemma}
\label{lemmarearrange}
Given $u_{n_1}^1 u_{n_2}^2 \cdots u_{n_k}^k 1 \in V^{(s)}_N$, and a permutation $\sigma \in S_k$, we have:
\[ u_{n_1}^1 u_{n_2}^2 \cdots u_{n_k}^k 1 = u_{n_{\sigma(1)}}^{\sigma(1)} u_{n_{\sigma(2)}}^{\sigma(2)} \cdots u_{n_{\sigma(k)}}^{\sigma(k)} 1  + R,\]
\noindent
for some $R \in V^{(s-1)}_N$.
\end{lemma}

\begin{proof}
It suffices to show this for transpositions. Commuting the $i$-th and $i+1$-st modes we have:
\begin{eqnarray*}
u_{n_1}^1 u_{n_2}^2 \cdots u_{n_{i-1}}^{i-1}u_{n_{i+1}}^{i+1} u_{n_{i}}^{i}u_{n_{i+2}}^{i+2} \cdots  u_{n_k}^k \1 =\\
u_{n_1}^1 u_{n_2}^2 \cdots u_{n_i}^i u_{n_{i+1}}^{i+1} \cdots  u_{n_k}^k \1 \\
- u_{n_1}^1 u_{n_2}^2 \cdots u_{n_{i-1}}^{i-1} [u_{n_i}^i, u_{n_{i+1}}^{i+1}] u_{n_{i+2}}^{i+2} \cdots  u_{n_k}^k \1.
\end{eqnarray*}
The two terms without the commutator are in the same filtration level as they have the same length and same modes. If the original string is in the $s$-th filtration level, then
\begin{equation*}
\sum_{t \ge 0}^k \wt{u^{t}}  + (N-1)k \le s
\end{equation*}
where the $i$-th and $i+1$-st modes contribute $\wt u^i +\wt u^{i+1} +2N-2$ to the filtration level $s$. By \eqref{lemmacommutator}
\begin{equation*}
[u_{n_i}^i, u_{n_{i+1}}^{i+1}]=\sum_{j \ge 0} \begin{pmatrix}n_i\\j\end{pmatrix} (u^{i}_ju^{i+1})_{n_i+n_{i+1}-j}.
\end{equation*}
The $j$-th term in this sum contributes $\wt{(u_j^i u^{i+1})}+N-1$ to the filtration level. So the filtration level of $u_{n_1}^1 u_{n_2}^2 \cdots u_{n_{i-1}}^{i-1} [u_{n_i}^i, u_{n_{i+1}}^{i+1}] u_{n_{i+2}}^{i+2} \cdots  u_{n_k}^k \1$ is
\[ \wt{u^i_ju^{i+1}} + \sum_{t \ge 0, t \ne i, i+1}^k \wt{u^{t}} + (N-1)(k-1)\]
\noindent
where
\begin{eqnarray} \label{weight}
\wt{(u_j^i u^{i+1})}&=&\wt{u^i} + \wt{u^{i+1}}+N-j-2\\
&<& \wt{u^i} + \wt{u^{i+1}} +2N-2 \nonumber
\end{eqnarray}
for $N \ge1$. So 
\[\wt{u^i_ju^{i+1}} + \sum_{t \ge 0, t \ne i, i+1}^k \wt{u^{t}} + (N-1)(k-1)<s
\]
and $u_{n_1}^1 u_{n_2}^2 \cdots u_{n_{i-1}}^{i-1} [u_{n_i}^i, u_{n_{i+1}}^{i+1}] u_{n_{i+2}}^{i+2} \cdots  u_{n_k}^k \1 \in V^{(s-1)}_N$ .
\end{proof}

\begin{remark}
In Equation (\ref{weight}), $\wt{(u_j^i u^{i+1})}<\wt{u^i} + \wt{u^{i+1}}$ for any $N \ge 1$.  In particular, the remainder term $R$ in Lemma \ref{lemmarearrange}  is in $V^{(s-1)}_{M}$ for any $M \ge 1$.
\end{remark}

A simple application of the above result is used in the proof of the main theorem. We can rewrite any monomial containing a mode with a nonnegative index in terms of monomials of a strictly lower filtration level.

\begin{lemma}
\label{ignorepositiveslemma}
Let $u^i \in V$ and $n_i \in \Z$ for $i = 1, \cdots k$, such that $u_{n_1}^1 u_{n_2}^2 \cdots u_{n_k}^k 1 \in V^{(s)}_N$, with $\sum_{i=1}^k \wt{u^i} + k(r-1) = s$. If $n_j \ge 0$ for some $j  \le k$, then $u_{n_1}^1 u_{n_2}^2 \cdots u_{n_k}^k \1 \in V^{(s-1)}_N$.
\end{lemma}

\begin{proof}
In Lemma \ref{lemmarearrange} above, set the permutation $\sigma \in S_k$ to be such that $\sigma(j) = k$. The result then follows from the creation axiom. 
\end{proof}

The next lemma shows that difference-$N$ filtrations satisfy a replacement property.  We may replace $u^i$ in a monomial $u_{n_1}^1 u_{n_2}^2 \cdots u_{n_k}^k \1$ with any representative of the coset $u^i+C_{N+1}(V)$ of $V/C_{N+1}(V)$, and the result is a monomial equivalent to the original modulo terms from lower filtration levels.  We will use this result to restrict the number of generators of the difference-$N$ spanning set.

\begin{lemma}
\label{lemmareplace}
Given $u_{n_1}^1 u_{n_2}^2 \cdots u_{n_k}^k \1 \in V^{(s)}_N$, 
\[ u_{n_1}^1 u_{n_2}^2 \cdots u_{n_k}^k \1 =  x_{n_1}^1 x_{n_2}^2 \cdots x_{n_k}^k \1  + R,\]
\noindent
where $x^i$ is a representative of $u^i+C_{N+1}(V)$ and $R \in V^{(s-1)}_N$.
\end{lemma}

\begin{proof}
We start with $u^i = x^i  + \sum a_{-N-1}b$. Linearity gives:
\begin{eqnarray*}
u_{n_1}^1 u_{n_2}^2 \cdots u_{n_i}^i \cdots u_{n_k}^k 1 &=&u_{n_1}^1 u_{n_2}^2 \cdots x_{n_i}^i \cdots u_{n_k}^k \1\\
&+& \sum u_{n_1}^1 u_{n_2}^2 \cdots (a_{-N-1} b)_{n_i} \cdots u_{n_k}^k \1
\end{eqnarray*}
\noindent
with $\wt{a_{-N-1}b} = \deg{a_{-N-1}} + \wt{b} = \wt{a} + \wt{b} +N$. However, using the associativity identity \eqref{lemmaassociator}, we can rewrite $a_{-N-1}b$ as:
\begin{eqnarray*} (a_{-N-1} b)_{n_i} &=& \sum_{r \ge 0} \begin{pmatrix} {-N-1} \\ r \end{pmatrix} (-1)^r a_{{-N-1}-r} b_{{n_i}+r} \\
&-&  \sum_{r \ge 0} \begin{pmatrix} {-N-1} \\ r \end{pmatrix} (-1)^{{-2}-r}b_{{-_N-1}+{n_i}-r} a_r.
\end{eqnarray*} 
\noindent
The contribution to the difference-$N$ filtration level of terms of the form $a_{-(N+1)}b$ is:
\begin{eqnarray*}
\wt{a_{-(N+1)}b}&=&\wt{a} + \wt{b} + ((N+1)-1) + (N-1)\\
&=&\wt{a} + \wt{b} + 2N-1.
\end{eqnarray*}
The contribution to the  filtration level of terms of the form $a_{p}b_{q}$ is: 
\[
(\wt{a}+(N-1)) + (\wt{b} +(N-1)) = \wt{a} + \wt{b} +2N-2.
\]
Therefore, terms of the form $\sum u_{n_1}^1 u_{n_2}^2 \cdots (a_{-N-1} b)_{n_i} \cdots u_{n_k}^k \1$ are in $V^{(s-1)}_N$.
\end{proof}

\begin{remark}
The remainder term $R$ in Lemma \ref{lemmareplace} is in $V^{(s-1)}_M$ for all $M \le N$ since the difference-$N$ filtration is finer than difference-$M$ filtration for $M \le N$. 
\end{remark}

In summary, the commutativity identity for modes of a M\"{o}bius vertex algebra ensures that the difference-$N$ filtrations satisfy a rearrangement property.  The associativity identity ensures that these filtrations satisfy a replacement property.    The associativity identity proves to be even more useful, as it is used in the next section to construct straightening identities. 

\section{Straightening identities}

In this section we develop the identities used to impose the desired difference conditions on spanning set elements.  We call these identities \emph{straightening identities}, using the terminology of Lepowsky and Wilson in \cite{MR752821}.  The straightening identities of this section only require application of the associativity identity, and so apply in general to vertex algebras.  However some of the analysis in remarks following the straightening lemmas will assume that elements $u$ and $v$ are graded.

\subsection{Difference-one straightening identity} 

This is Lemma 7 from \cite{MR1990879}, and is used to replace any spanning set expression involving repeated indices of modes with expressions without repeated indices. In this section we generalize this identity and its associated properties.

\begin{lemma}[{\bf Difference-one straightening identity}]
\label{lemmastraight1}
For $u,v \in V$, a vertex algebra, and $n \in \Z$ with $n \le 0$,
\[ u_{n}v_{n} = (u_{-1}v)_{2n+1} - \sum_{i \ge 0, i \neq -n} u_{-i} v_{2n+i} 
- \sum_{i \ge 0} v_{2n-i} u_i. \]
\end{lemma}

\begin{proof}
Compute $(u_{-1}v)_{2n+1}$ using the associativity identity \eqref{lemmaassociator}, and isolate $u_{n}v_{n}$.
\end{proof}

Germane properties of this identity are: substituting the right-hand expression for the the left-hand one in a monomial preserves both the filtration level and weight of the monomial, the modes on the left-hand side of the identity have no ``close" indices, and with the exception of the associator term for $n=1$, each mode on the left-hand side has a smaller index than either mode on the right-hand side of the identity.

\begin{remark}[{\bf Preserves Filtration Level}]
For a M\"{o}bius vertex algebra, all terms that appear in the above identity have the same contribution for the difference-one filtration level.  All the $u_j v_k$ terms contribute $\wt{u} + \wt{v}$ to the difference-one filtration level, and $(u_{-1}v)_{2n+1}$ contributes $\wt{u} +\wt{v}$ to the difference-one filtration level as well.  For the filtration for larger $N$, the product terms $u_iv_j$ contribute  $ \wt{u} + \wt{v} + 2(N-1)$ to the difference-N filtration level. The associator term $(u_{-1}v)_{2n+1}$ contributes $\wt{u} +\wt{v}+N-1$ to the difference-$N$ filtration level because it only has a length of 1.  So for the difference-one filtration, the terms on the right-hand side of the straightening identity contributes the same amount to the filtration level as the term on the left-hand side.  This holds for any difference-N filtration with the exception of the associator term, which contributes a strictly smaller amount to the filtration level for $N>1$.
\end{remark}

\begin{remark}[{\bf Preserves Degrees}]
Each of the products of modes and the associator mode that appear in the difference-one spanning set have the same degrees as operators on $V$.  As a result, substitution the left-hand side expression for the right-hand one in a monomial preserves the weight of that monomial.
\end{remark}

\begin{remark}[{\bf Smaller Indices}]
In the product terms $u_jv_k$ that appear on the right-hand side of the straightening identity above, either $j$ or $k$ is less than $n$.
\end{remark}

\begin{remark}[{\bf No Small Differences}]
The first sum in the difference-one straightening lemma has no product terms with indices of difference one or less, and the second sum has no product terms with difference one or less for $n<0$.
\end{remark}

The goal of the remainder of this section is to develop a family of identities that share similar features.

\subsection{Associativity identities}
\label{Subsection_NotationAB}
For the higher difference straightening identities, we use multiple instances of the associativity identity to create an identity that eliminates modes with indices that are ``too close".  For even $N$, we look at the associativity identities for $(u_{-1-r}v)_{2n+r}$ where $r$ ranges between 0 and $N-1$. For odd $N$, we look at the identities for $(u_{-1-r}v)_{2n+1+r}$, where, once again, $r$ ranges between 0 and $N-1$. Here we introduce some notation for the summation terms that appear on the right-hand side of these particular associativity identities.

For $N=2k$ we examine the associativity identities of the form
\begin{eqnarray*}
(u_{-1-r}v)_{2n+r} &=&\sum_{i \ge 0} \binom{i+r}{i} u_{-1-r-i} v_{2n+r+i} \\
&&-  (-1)^{-1-r} \sum_{i \ge 0} \binom{i+r}{i}  v_{2n-1-i} u_i.
\end{eqnarray*}
Isolating from the first sum the pairs of modes with indices of difference-$N$ or less and reindexing these terms, we obtain 
\begin{eqnarray} \label{evenassoc_first}
(u_{-1-r}v)_{2n+r} &=&\sum_{i} \binom{i+r}{i} u_{-1-r-i} v_{2n+r+i}  \nonumber \\
&&+ \sum_{j = 0}^{2k-1} \binom{-n+k-1-j}{-n+k-1-j-r} u_{n-k+j} v_{n+k-1-j} \nonumber \\
&&-  (-1)^{-1-r} \sum_{i \ge 0} \binom{i+r}{i} v_{2n-1-i} u_i,
\end{eqnarray}
where in the first sum $i \in \N \backslash \{-n-k-r, \cdots, -n+k-r-1\}$. We call this sum $A_{r,n,k}(u,v)$. In other words, we define:
\begin{eqnarray*}
A_{r,n,k}(u,v)&=&\sum_{i\ge 0} \binom{i+r}{i} u_{-1-r-i} v_{2n+r+i} \\
&&- \sum_{j = 0}^{2k-1} \binom{-n+k-1-j}{-n+k-1-j-r} u_{n-k+j} v_{n+k-1-j}.  
\end{eqnarray*}
We call the second sum in \eqref{evenassoc_first} $B_{r,n}(u,v)$. In other words we define: 
\[B_{r,n}(u,v)=(-1)^{-r}\sum_{i\ge 0} \binom{i+r}{i} v_{2n-1-i} u_{i}.\]
\noindent
After rearranging and using this notation, we rewrite (\ref{evenassoc_first}) as  
\begin{eqnarray} \label{evenassoc_withAB}
\sum_{j = 0}^{2k-1} \binom{-n+k-1-j}{-n+k-1-j-r} u_{n-k+j} v_{n+k-1-j} \nonumber \\
=(u_{-1-r}v)_{2n+r} - A_{r,n,k}(u,v) - B_{r,n}(u,v).
\end{eqnarray}

For $N=2k+1$ we rewrite the associativity identity for $(u_{-1-r}v)_{2n+1+r}$ in a similar way to obtain
\begin{eqnarray}\label{oddassoc_withAB}
\sum_{j = 0}^{2k} \binom{-n+k-1-j}{-n+k-1-j-r} u_{n-k+j} v_{n+k-j} \nonumber\\
=(u_{-1-r}v)_{2n+1+r} - A'_{r,n,k}(u,v) - B'_{r,n}(u,v)
\end{eqnarray}
where 
\begin{eqnarray*} 
A'_{r,n,k}(u,v)&=&\sum_{i\ge 0} \binom{i+r}{i} u_{-1-r-i} v_{2n+1+r+i}  \\
&&- \sum_{j = 0}^{2k} \binom{-n+k-1-j}{-n+k-1-j-r} u_{n-k+j} v_{n+k-j}  
\end{eqnarray*}
and 
\[
B'_{r,n}(u,v)=(-1)^{-r}\sum_{i\ge 0} \binom{i+r}{i} v_{2n-i} u_{i}.
\]

\begin{example}
Using this new notation we may rewrite the difference-one straightening identity (Lemma \ref{lemmastraight1}) as
\[
u_{n}v_{n} = (u_{-1}v)_{2n+1} - A'_{0,n,0}(u,v) - B'_{0,n}(u,v).
\]
\end{example}

Below we make some observations about these rearranged associativity identities \eqref{evenassoc_withAB} and \eqref{oddassoc_withAB} analogous to our earlier observations about the difference-one straightening idenitity.

\begin{remark}[{\bf Preserves Filtration Level}]\label{UsefulRemark1}
Each of the product terms that appears in a given $A$, $A'$, $B$ or $B'$ 
expression contributes the same amount to the filtration level of a monomial, for any difference-N filtration chosen.  This is true because the length and vectors are unchanged, just rearranged with difference indices of modes. 

In each version of the associativity identity written in the form of \eqref{evenassoc_withAB} or \eqref{oddassoc_withAB}, the associator term $u_{-1-r}v$ contributes $\wt u  - \wt v + r - (N-1)$ to the difference $N$ filtration level.  However it is shorter so when it replaces a product term it effectively changes the filtration level by $\wt u  - \wt v + r - 2(N-1)$.  If $r \le N-1$ then $u_{-1-r}v$ contributes the same or a lesser amount to the filtration level than $u_iv_j$.
\end{remark}

\begin{remark}[{\bf Preserves Degrees}]\label{UsefulRemark2}
Each of the products of modes and the associator mode that appear in \eqref{evenassoc_withAB} and \eqref{oddassoc_withAB} have the same degrees as operators on $V$, as is true for the associativity identity in general.
\end{remark}

\begin{remark}[{\bf Smaller Indices}]\label{diffnstraightindex}
For $N =2k$, in each of the product terms $u_iv_j$ that appear within the sums $A_{0,n,k}$, $A_{1,n,k}$ $\cdots$ $A_{2k-1,n,k}$,
either $i$ or $j$ is less than $n-k$. Similarly, for $N=2k+1$, in each of the product terms $u_iv_j$ that appear within the sums $A^{\prime}_{0,n,k}$, $A^{\prime}_{1,n,k}$ $\cdots$ $A^{\prime}_{2k,n,k}$,
either $i$ or $j$ is less than $n-k$. 
\end{remark}

\begin{remark}[{\bf No Small Differences}]\label{UsefulRemark4}
$A_{r,n,k}(u,v)$ has no product term with indices of difference $2k$ or less, and $B_{r,n}(u,v)$ has no product terms with difference $2k$ or less for $-2n+1<2k$ or $-n \le k-1$.  $A'_{r,n,k}(u,v)$ has no product term with indices of difference $2k$ or less, and $B_{r,n}(u,v)$ has no product terms with difference $2k+1$ or less for $-2n<2k+1$ or $-n \le k$. 
\end{remark}

\subsection{Arbitrary Difference Straightening Identities}

This most general case involves $N$ instances of the associativity identity \eqref{lemmaassociator}, each of which gives us a linear equation. The complete system of $N$ equations has a nonsingular coefficient matrix, allowing us to solve for any term whose modes are less than $N$ apart. We consider the even and odd cases separately.

\begin{lemma}[{\bf Difference-$N$ straightening identity}]
\label{lemmastraightn}
Let $V$ be a vertex algebra and $u,v \in V$. Let $n \in \Z_-$.  If $N = 2k$ for some $k \in \Z_+$, such that $n+k-1<0$, then
\[
u_{n-k}v_{n+k-1} = \sum_{r=0}^{2k-1} c_N(r,n)\left((u_{-1-r}v)_{2n+r}- A^{}_{r,n,k}(u,v) - B_{r,n}(u,v)\right)\]
\noindent
for some coefficients $c_N(r,n)$. If $N = 2k+1$ for some $k \in \Z_+$ and $n+k<0$, then
\[
u_{n-k}v_{n+k} = \sum_{r=0}^{2k} c_N(r,n)\left((u_{-1-r}v)_{2n+1+r} - A^{\prime}_{r,n,k}(u,v) - B^{\prime}_{r,n}(u,v) \right)\]
\noindent
for the same coefficients $c_N(r,n)$. 
\end{lemma}

\begin{proof}
We prove the lemma for $N$ even; the proof for odd $N$ is similar.  Recall Equation \eqref{evenassoc_withAB} for $0 \le r \le N-1$:
\begin{eqnarray*} 
\sum_{j = 0}^{2k-1} \binom{-n+k-1-j}{-n+k-1-j-r} u_{n-k+j} v_{n+k-1-j} \\
=(u_{-1-r}v)_{2n+r} - A_{r,n,k}(u,v) - B_{r,n}(u,v).
\end{eqnarray*}
\noindent 
Each of these $N$ equations is linear in the product terms $u_{n-k+j} v_{n+k-1-j}$.  The coefficient matrix of the corresponding linear system is the $N \times N$ matrix $S_N(-n+k-1)$ whose $ij$-th entry is $\binom{-n+k-1-(j-1)}{-n+k-1-(i-1)-(j-1)}$.  For $n+k-1<0$, this matrix is invertible; see Proposition \ref{LUfactor} in the appendix.  

Solving for the product term $u_{n-k} v_{n+k-1}$ now yields
\[
u_{n-k} v_{n+k-1}= 
\sum_{r=0}^{2k-1} c_N(r,n) \left(\sum_{j = 0}^{2k-1} \binom{-n+k-1-j}{-n+k-1-j-r} u_{n-k+j} v_{n+k-1-j} \right)
\]
\noindent
where 
\[c_N(r,n)= -\sum_{t=1}^N (-1)^{r+t} \binom{-n+k-1-r}{-n+k-1-(t-1)}\]
\noindent
as detailed in the appendix.  Finally, substituting the right-hand side of Equation \eqref{evenassoc_withAB} into this equation  
gives the desired result.
\end{proof}

\begin{remark}
In the product term making up the left-hand side of the difference-$N$ straightening identity, the difference between the two modes is $N-1$.  For each product term on the of the right-hand side of the identity, the indices of the modes have difference $N$ or more.  So one applies a difference-$N$ straightening identity to two modes of difference $N-1$ to rewrite it in terms of modes that satisfy a larger difference condition.
\end{remark} 

\begin{remark}
We restricted our lemma to the case $n+k-1<0$ for $N=2k$ and to the case $n+k<0$ for $N=2k+1$. We do not worry about the remaining cases ($n+k-1\ge 0$ and $n+k\ge 0$, respectively) because of Lemma \ref{ignorepositiveslemma}.  In other words, we will need to apply the difference-$N$ straightening identity only to products of modes of the form $u_iv_j$ with $j<0$. 
\end{remark}

The following remarks generalize the remarks following the difference-one straightening identity.

\begin{remark} [{\bf Preserves Filtration Level}] \label{straightrem1}
Following Remark \ref{UsefulRemark1} substituting the right-hand side of a difference-$N$ straightening identity for the left will not introduce any terms of higher filtration for any difference-$M$ filtration with $M \ge N$. 
\end{remark}

\begin{remark}[{\bf Preserves Degrees}]
Following Remark \ref{UsefulRemark2} all associator terms and product terms have the same degrees as operators on a M\"{o}bius vertex algebra $V$. Thus in any given monomial, substituting the right-hand side of a difference-$N$ straightening identity for the left preserves the weight of the monomial.
\end{remark}

\begin{remark}[{\bf Smaller Indices}]
Following Remark \ref{diffnstraightindex} each product term that appears on the right-hand side of the ``even'' straightening identity for $u_{n-k} v_{n+k-1}$ has at least one mode whose index is smaller than $n-k$.  This is true for the ``odd'' straightening identities as well.
\end{remark}

\begin{remark}[{\bf No Small Differences}] Following Remark \ref{UsefulRemark4} each product term on the left-hand side of a difference-$N$ straightening identity consists merely of modes whose indices differ by $N$ or more.
\end{remark}

\begin{example}
To get a straightening identity for $N = 2$ ($N = 2k$ with $k=1$), we need to use two rearrangements of the associativity identity, namely those for $(u_{-1}v)_{2n}$, corresponding to $r=0$, and $(u_{-2}v)_{2n+1}$, corresponding to $r=1$:
\[
u_{n-1}v_{n} + u_{n}v_{n-1} =
(u_{-1}v)_{2n} - A_{0,n,1}(u,v) - B_{0,n}(u,v),\]
\noindent
and:
\[
(-n) u_{n-1}v_{n} + (-n-1) u_{n}v_{n-1}  =
(u_{-2}v)_{2n+1} - A_{1,n,1}(u,v) + B_{1,n}(u,v). \]
\noindent 
Solving for $u_{n-1}v_n$ yields the difference-two straightening identity
\begin{eqnarray*}
u_{n-1}v_n &=& (n+1)((u_{-1}v)_{2n} - A_{0,n,1}(u,v) - B_{0,n}(u,v))  \\
&+& ((u_{-2}v)_{2n+1} - A_{1,n,1}(u,v) + B_{1,n}(u,v))
\end{eqnarray*}
which is valid for $u,v \in V$ and $n \in \Z$ with $n \le -1$.
 
\end{example}

\begin{example}
Similarly, the identity for $N = 3$ ($N = 2k+1$ with $k=1$) is constructed from the associativity identities for $(u_{-1}v)_{2n+1}$, $(u_{-2}v)_{2n+2}$, and $(u_{-3}v)_{2n+3}$:
\begin{eqnarray*}
u_{n-1}v_{n+1} + u_{n}v_{n}+ u_{n+1}v_{n-1} &=& \\
(u_{-1}v)_{2n+1} - A'_{0,n,1}(u,v) - B'_{0,n}(u,v), \\
\\
(-n)u_{n-1}v_{n+1} + (-n-1)u_{n}v_{n}+ (-n-2)u_{n+1}v_{n-1} &=& \\
(u_{-2}v)_{2n+2} - A'_{1,n,1}(u,v) - B'_{1,n}(u,v), \\
\\
\tfrac{n^2+n}{2}u_{n-1}v_{n+1} + \tfrac{n^2+3n+2}{2}u_{n}v_{n}+ \tfrac{n^2+5n+6}{2}u_{n+1}v_{n-1} &=& \\
(u_{-3}v)_{2n+3}- A'_{2,n,1}(u,v) - B'_{2,n}(u,v).
\end{eqnarray*}
 \noindent
Solving for $u_{n-1}v_{n+1}$ yields the difference-3 straightening identity
\begin{eqnarray*}
u_{n-1}v_{n+1} &=& \tfrac{n^2+3n+2}{2}((u_{-1}v)_{2n+1} - A'_{0,n,1}(u,v) - B'_{0,n}(u,v))\\
&&+ (n+2)((u_{-2}v)_{2n+2} - A'_{1,n,1}(u,v) - B'_{1,n}(u,v))\\
&&+ ((u_{-3}v)_{2n+3}- A'_{2,n,1}(u,v) - B'_{2,n}(u,v))
\end{eqnarray*}
which is valid for $u,v \in V$ and $n \in \Z$ with $n \le -2$.
\end{example}

\section{Spanning sets of arbitrary difference}
\label{SectionMainTheorem}

Using the straightening lemmas from the previous section, we will now construct Poincar\'{e}-Birkhoff-Witt-like spanning sets for M\"{o}bius vertex algebras that satisfy a difference-$N$ condition for any $N >0$.  That is, we will show that every monomial of finite length can be rewritten as a sum of finitely many monomials, each of which satisfies the desired ordering restriction.  The proof of this result is a generalization of the difference-one spanning set for twisted modules for vertex operator algebras \cite{MR2039213}. 

The following theorem, the difference-$N$ spanning set theorem for M\"{o}bius vertex algebras, is the main result of this paper:

\begin{theorem} \label{diffnthm}
For a $\N$-graded M\"{o}bius vertex algebra $V$ and $N\in\Z_+$, let $X$ be a set of homogeneous representatives of a spanning set for $V/C_{N+1}(V)$,  $V$ is spanned by monomials of the form
\begin{eqnarray}
x^1_{n_1} x^2_{n_2} \cdots x^k_{n_k}\1
\end{eqnarray}
where $k \in \N$; $x^{1}, \ldots, x^{k} \in X$; $n_1, \ldots , n_k \in \Z_{-}$; and $n_1 < n_2 < \cdots < n_k < 0$ with $n_i - n_{i+1} \ge N$ for each $1\ge i \ge k-1$.
\end{theorem}

The proof uses induction and a recursive algorithm for rewriting monomials in the desired form.  We will separate the proof that the recursive algorithm terminates after a finite number of iterations (Lemma \ref{terminate}) from the proof of the rest of Theorem \ref{diffnthm}. 

\begin{proof}[Proof of the rest of Theorem \ref{diffnthm}]
The proof is based on an induction argument on pairs $(s,k) \in \N \times \N$ with the order
\begin{eqnarray*}
(s,k) < (s',k') \iff s < s' \textmd{ or } s = s' \textmd{ and } k < k',
\end{eqnarray*}
where $s$ is the filtration level and $k$ is the length of a given monomial of the form  $u^1_{m_1} u^2_{m_2} \cdots u^k_{m_k}\1$.

For fixed, arbitrary $N$, let $X_N$ be a set of representatives of a basis for $V/C_{N+1}(V)$. The $(s,k)$-th induction hypothesis for $N$ is:

{\it For any $(s,k) < (s',k')$, every element of $V^{(s)}_N$ of length $k$ can be rewritten as a linear combination of  monomials in $V^{(s-1)}_N$ and monomials in $V^{(s)}_n$ of the form $x^1_{m_1} x^2_{m_2} \cdots x^{k}_{m_{k}}\1$ with $x^i \in X_N$ for all $i$, and $m_i < 0$ with $m_{i+1} - m_i \ge N$ for each $i$.}

We first show that the induction hypothesis is true for base cases:  $(0,0)$ and $(s,1)$. The monomials of length 0 are $c\1$ where $c$ is a scalar; these are already in the desired form.  Strings of length one have the form $u_n\1$. If $u_n\1\in V^{(s)}_N$, applying the replacement lemma (Lemma \ref{lemmareplace}) yields $x_n\1+R$ where $R \in V^{(s-1)}_N$, satisfying the induction hypothesis.  

Assume now that the $(s,k)$-th induction hypothesis is true and consider a generic string of filtration level $s$ and length $k$: $u^1_{n_1} u^2_{n_2} \cdots u^{k}_{n_{k}}1$.  Applying Lemma \ref{lemmareplace} replaces each of the modes in the monomial with representatives from $X_n$. Applying the reordering lemma (Lemma \ref{lemmarearrange}) places the modes in increasing order, yielding $x^1_{m_1} x^2_{m_2} \cdots x^{k}_{m_{k}}\1 + R$ where the first string is in $V^{(s)}_N$ and $R \in V^{(s-1)}_N$.  We may apply the induction hypothesis to $R$. Therefore, it is sufficient to restrict attention to monomials $x^1_{n_1} x^2_{n_2} \cdots x^{k}_{n_{k}}\1 \in V^{(s)}_N$ where the $x^i \in X_N$ and $n_1 \le n_2 \le \cdots \le n_{k} < 0$.   By Lemma \ref{ignorepositiveslemma}, if $n_i \ge 0$ for some $i$, then the monomial is in $V^{(s-1)}_N$.

Applying the induction hypothesis once again to the $(k-1)$-length tail of the above string yields a monomial of the form 
\begin{equation} \label{firstmonomial}
x^1_{m_1} x^2_{m_2} \cdots x^{k}_{m_{k}}\1,
\end{equation}
where $x^1_{m_1} x^2_{m_2} \cdots x^{k}_{m_{k}}\1 \in V^{(s)}_N$, $x^i \in X_N$,  $m_i < 0$ with $m_{i+1} - m_i \ge N$ for each $i \ge 2$, and $m_1=n_1$, plus an element in a strictly lower filtration level. We may then apply the induction hypothesis to this lower filtration element and further restrict our attention to monomials whose modes satisfy the difference-N condition with the possible exception of the first mode. 

At this point, there are three cases for $m_1$ and $m_2$:
\begin{enumerate}
\item $m_2 - m_1 \ge N$;
\item $0 \le m_2-m_1 < N$; and
\item $m_2-m_1<0$.
\end{enumerate}

\noindent\textbf{Case 1:}\\
In the first case, $m_2 - m_1 \ge  N$, the string $x^1_{m_1} x^2_{m_2} \cdots x^{k}_{m_{k}}\1$ is already in the desired form, with all the modes satisfying the difference-$N$ condition, and we are done.

\noindent\textbf{Case 2:}\\
In the second case, $0 \le m_2 - m_1 < N$, we apply the difference-$M$ straightening lemma where $M=m_2 - m_1+1$ to the pair $x^1_{m_1}x^2_{m_2}$. This allows us to rewrite the whole monomial  $x^1_{m_1} x^2_{m_2} \cdots x^{k}_{m_{k}}\1$ in terms of monomials which are composed of some monomials of equal length and others which are shorter. In the following analysis, we assume that $M$ is even. The analysis for odd $M$ is similar.

Applying the difference-$M$ straightening identity to the monomial yields:
\begin{eqnarray}
x^1_{m_1}x^2_{m_2} \cdots x^{k}_{m_{k}}\1 
&=& \sum_{r=0}^{2k-1} c_N(r,n)(x^1_{-1-r}x^2)_{2n+r}x^3_{m_3} \cdots x^{k}_{m_{k}}\1 \\
&&- \sum_{r=0}^{2k-1} c_N(r,n)A_{r,n,k}(x^1,x^2)x^3_{m_3} \cdots x^{k}_{m_{k}}\1 \\
&&- \sum_{r=0}^{2k-1} c_N(r,n)B_{r,n}(x^1,x^2)x^3_{m_3} \cdots x^{k}_{m_{k}}\1
\end{eqnarray}
where $n=\frac{m_1+m_2+1}{2}$ and $k=\frac{m_2 -m_1+1}{2}$.

In the first sum, each monomial has an associator term at the beginning and is shorter than the original monomial. By Remark \ref{straightrem1} each term in the first sum has the same or lower filtration level.  Because the shorter strings lie in, at most, the same filtration level as \eqref{firstmonomial}, we apply the induction hypothesis to rewrite these terms in the desired form.

This leaves the strings of length $k$ in the second and third sums which have the form $y^1_{n_1} y^2_{n_2} x^3_{m_3} \cdots x^{k}_{m_{k}}\1$ with $|n_2 - n_1| \ge M$ and $m_3 < m_4 < \cdots < m_{k} < 0$ such that $m_{i+1} - m_i \ge n$ for each $i\ge 3$.  Recall that each product in the sums involved in an application of a straightening lemma contains at least one mode which is less than the initial mode of the original product, see Remark \ref{diffnstraightindex}. This implies that at least one mode, either $y^1_{n_1}$ or $y^2_{n_2}$, has index less than $m_1$. Because of the rearrangement lemma (Lemma \ref{lemmarearrange}), we may assume without loss of generality that $n_1$ is less than $m_1$.

Now we do not know the relationship between $n_2$ and $m_3$.  We  apply the induction hypothesis to the tail $y^2_{n_2} x^3_{m_3} \cdots x^{k}_{m_{k}}1$ to rewrite our string in terms of new monomials whose tails satisfy the difference-$N$ condition. This process produces a new monomial of the same form as \eqref{firstmonomial},  but with $n_1<m_1$.  That is, after relabeling this process results in a new monomial
\begin{equation} \label{secondmonomial}
x^1_{m'_1} x^2_{m'_2} \cdots x^{k}_{m'_{k}}\1
\end{equation}
which lies in $ V^{(s)}_n$ where the $x^i$ are in the set $X_n$ of distinguished elements, $m'_1 = n_1$, and $m'_2 < m'_3 < \cdots < m'_{k} < 0$ with $m'_{i+1} - m'_i \ge n$ for each $i> 1$. However the initial mode $m'_1$ of this second string is less than the initial mode $m_1$ of the original \eqref{firstmonomial}.  So we return to the scenario with the three cases for the relationships between the indices of first two modes, but the index of the first mode in this new monomial is strictly smaller. Note that the new monomial has the same weight, length, and filtration as \eqref{firstmonomial}.

\noindent\textbf{Case 3:}\\
In the third case, $m_2-m_1<0$, or equivalently, $m_2 < m_1$. We apply Lemma \ref{lemmarearrange} to our string and end up with $x^2_{m_2} x^1_{m_1} x^3_{m_3} \cdots x^{k}_{m_{k}}1$ plus an element of a lower filtration level.  We may rewrite the element of the lower filtration level in the desired form using the induction hypothesis.  We then apply the induction hypothesis to the $(k-1)$-length tail of $x^2_{m_2} x^1_{m_1} x^3_{m_3} \cdots x^{k}_{m_{k}}1$.  This results in a term in a lower filtration level and a monomial of the form $y^1_{n_2} y^2_{n_1}  \cdots y^{k}_{n_{k}}1$, with $y^1_{n_2}=x^2_{m_2}$, whose tail satisfies the difference-N condition.  Again the induction hypothesis can be applied to the lower filtration term. This once again returns us to the situation where we have a string and there are three possible cases for the relationship between the first two modes. However, now the index of the first mode is smaller than that of the original monomial. We note that the new monomial has the same weight, length, and filtration as \eqref{firstmonomial}.

In each of the three cases, we rewrite our initial string in terms of strings of the desired form by applying the induction hypothesis and end up with terms of same length, same weight, same filtration level, and first mode with a strictly smaller index.  Lemma \ref{terminate} ensures that this process of rewriting strings using this process eventually terminals, leaving only sums of monomials satisfying the difference-$N$ condition. Thus, modulo the proof of Lemma \ref{terminate}, we are done.
\end{proof}

In order to complete the proof of Theorem \ref{diffnthm}, we now show that the process of rewriting monomials eventually terminates. Here is the precise statement of this fact:

\begin{lemma} \label{terminate}
The process in the proof of Theorem \ref{diffnthm} of rewriting monomials of the form
\[x^1_{m_1} x^2_{m_2} \cdots x^{k}_{m_{k}}\1\]
using Lemmas \ref{lemmarearrange}, \ref{lemmareplace}, and \ref{lemmastraightn} terminates after a finite number of iterations.
\end{lemma}

\begin{proof}
Both the original monomial of the form \eqref{firstmonomial} and any new monomial arising from an application of either of the Lemmas \ref{lemmarearrange}, \ref{lemmareplace}, and \ref{lemmastraightn} will have the same weight $w$, the same filtration level $s$, and same length $k$. More specifically, the weight of the monomial \eqref{firstmonomial} is 
\[\sum_{i=1}^k (\wt{x^i}-m_i-1)= s-Nk-\sum_{i=1}^k m_i.
\]
Similarly the weight of the resultant monomial of the form \eqref{secondmonomial} is 
\[s-Nk-\sum_{i=1}^k m'_i.
\]
Our earlier assertion then implies then that $\sum_{i=1}^k m'_i=\sum_{i=1}^k m_i$, 
where $\sum_{i=1}^k m_i$ is fixed by the original monomial \eqref{firstmonomial}.  If a process decreases $m'_1$ the remaining $m'_i$'s must increase to ensure that the sum is fixed. So for sufficiently small $m'_1$, $m'_k$ must be non-negative, and so $x^k_{m'_k}\1=0$.  This proves then that the process must terminate. \end{proof}

Theorem \ref{diffnthm} provides an alternate, direct proof that $C_N$-cofinite, $\N$-graded M\"{o}bius vertex algebras are finitely generated.  This was originally proved for certain vertex operator algebras in \cite{MR1990879}. It is now:

\begin{corollary} \label{diffncor}
Any $C_N$-cofinite, $\N$-graded M\"{o}bius vertex algebra is finitely generated for $N \ge 2$.
\end{corollary}

\begin{proof}
Let $V$ be a M\"{o}bius vertex algebra. If $V$ is $C_{N+1}$-cofinite, the set of generators $X_N$ of Theorem \ref{diffnthm} are finite.
\end{proof}

\section{Appendix}

In this section we establish some combinatorial identities arising in the construction of the difference-$N$ straightening lemma (Lemma \ref{lemmastraightn}).  In particular, we show that the coefficient matrix for a particular linear system is invertible and solve an instance of this system.

\begin{definition}
The {\it difference-$N$ straightening identity matrix for m} is the $N \times N$ matrix whose $ij$-th entry is $\binom{m-(j-1)}{m-(i-1)-(j-1)}$.  We denote this matrix as $S_N(m)$, where
\begin{eqnarray*}
S_N(m)=
\left[\begin{array}{cccc}
\binom{m}{m} & \binom{m-1}{m-1} & \cdots & \binom{m-(N-1)}{m-(N-1)} \\
\binom{m}{m-1} & \binom{m-1}{m-2} & \cdots & \binom{m-(N-1)}{m-N} \\
\vdots & \vdots & \ddots & \vdots \\
\binom{m}{m-(N-1)} & \binom{m-1}{m-N} & \cdots & \binom{m-(N-1)}{m-2(N-1)}
\end{array}\right].
\end{eqnarray*}
\end{definition}

\begin{definition}
The {\it upper triangular Pascal matrix} is the $N \times N$ matrix $P_N$ whose $ij$-th entry is $\binom{j-1}{i-1}$ where $1 \leq i,j \leq N$.  
\end{definition}

To prove the difference-$N$ straightening lemmas, it is necessary to prove that $S_N(m)$ non-singular for certain $m$.  We show that $P_N$ and $S_N(m)$ are row equivalent for $m \geq N-1$.  To prove this, we need:

\begin{lemma}\label{binomid}
For  $k,m \in \N$ with $k\ge 1$ and $m \geq k-1$,
\[\sum_{i=1}^k (-1)^{i-1}\binom{m-(i-1)}{m-(k-1)}\binom{m-(j-1)}{m-(i-1)-(j-1)}=\binom{j-1}{k-1}.\]
\end{lemma}

\begin{proof} We have:
\begin{eqnarray*}
&&\sum_{i=1}^k (-1)^{i-1}\binom{m-(i-1)}{m-(k-1)}\binom{m-(j-1)}{m-(i-1)-(j-1)} \\
&=& \sum_{i=1}^k (-1)^{i-1}\binom{m-(i-1)}{m-(k-1)}\binom{m-(j-1)}{(i-1)} \\
&=& (-1)^{2m-(k-1)}\binom{(k-1)-(j-1)-1}{k-1}
\end{eqnarray*}
by binomial identity (5.25) of \cite{MR1397498}:
\[\sum_{i \leq a} (-1)^{i}\binom{a-i}{b}\binom{c}{i-d}=(-1)^{a+b}\binom{c-b-1}{a-b-d}.
\]
Applying the upper negation identity 
\[\binom{a}{b}=(-1)^b\binom{b-a-1}{b},
\]
we get
\[(-1)^{k-1}\binom{(k-1)-(j-1)-1}{k-1}
=\binom{j-1}{k-1}.\]\end{proof}

\begin{proposition}
For $m \in \mathbb{N}$ and $m \geq N-1$, $S_N(m)$ is row equivalent to the upper triangular Pascal matrix $P_N$.
\end{proposition}

\begin{proof}
Let $\row_k(A)$ denote the $k$-th row of a matrix $A$.  We show that 
\begin{eqnarray*}
\row_k(S_N(m))+\sum_{i=1}^{k-1} (-1)^{i+k}\binom{m-(i-1)}{m-(k-1)}\row_i(S_N(m))\\=(-1)^{k+1}\row_k(P_N)
\end{eqnarray*}
or equivalently
\begin{equation}
\sum_{i=1}^{k} (-1)^{i-1}\binom{m-(i-1)}{m-(k-1)}\row_i(S_N(m))=\row_k(P_N) .\label{lincombo}
\end{equation}
In other words, the $k$-th row of $P_N$ is a linear combination of the first $k$ rows of $S_N(m)$.  Looking at the entries, we must show that
\[\sum_{i=1}^k (-1)^{i-1}\binom{m-(i-1)}{m-(k-1)}\binom{m-(j-1)}{m-(i-1)-(j-i)}=\binom{j-1}{k-1}
\]
which is precisely what Lemma \ref{binomid} states.
\end{proof}

\begin{proposition} \label{LUfactor}
For  $m \in \N$ and $m \geq N-1$, $\mbox{Det }S_N(m)=(-1)^{(N-1)}$.
\end{proposition}

\begin{proof}
Equation \eqref{lincombo} implies that
\begin{equation}\label{decomp}
L_N(m) S_N(m)=P_N
\end{equation}
where $L_N(m)$ is the $N \times N$ lower triangular matrix whose $ij$-th entry is $(-1)^{j-1}\binom{m-(j-1)}{m-(i-1)}$.  The determinant of $L_N(m)$ is $(-1)^{\frac{(N-1)N}{2}}$, and the determinant of $P_N$ is $1$, so the determinant of  $S_N(m)$ is $(-1)^{\frac{(N-1)N}{2}}$.
\end{proof}

In addition to showing $S_N(m)$ is non-singular for certain $m$, we can use Proposition \ref{LUfactor}, more specifically Equation \eqref{decomp},
to get the coefficients $c_N(r,n)$ of the $A$, $A'$, $B$, and $B'$ that appear in the difference-$N$ straightening identities of Lemma \ref{lemmastraightn}.  The coefficient $c_N(r,n)$ is the $r+1$-st component of the solution to the linear system $S_N(-n+k-1) {\bf x}= {\bf e_1}$: 
\[c_N(r,n)=\row_{r+1} S_N(-n+k-1)^{-1} {\bf e_1}
\]
 where is ${\bf e_1}$ is the first standard basis vector.  The matrix $P_N^{-1}$ is the upper triangular matrix whose $ij$-th entry is $(-1)^{i+j}\binom{j-1}{i-1}$. 
Therefore the $ij$-th component of $S_n(m)^{-1}$, which is $P_N^{-1}L_N(m)$, is:
\[\sum_{t=1}^N ((P_N)^{-1})_{it} (L_N(m))_{tj} = 
\sum_{t=i}^N (-1)^{i+j+t-1}\binom{t-1}{i-1}\binom{m-(j-1)}{m-(t-1)}.
\]
We only need the first row of this matrix. The $j$th entry in this first row is
\[\sum_{t=1}^N (-1)^{j+t} \binom{m-(j-1)}{m-(t-1)}.
\]
To find the coefficients in Lemma \ref{lemmastraightn}, we set $m = -n+k-1$ and $r=j-1$:
\[c_N(r,n) = -\sum_{t=1}^N (-1)^{r+t} \binom{-n+k-1-r}{-n+k-1-(t-1)}.
\]


\end{document}